\documentclass[12pt,a4paper]{article}
\usepackage{amssymb}

\usepackage{graphicx,amssymb,amsfonts,epsfig,amsthm,a4,amsmath,url}
\usepackage[latin1]{inputenc}

\newtheorem{ThmIntro}{Theorem}
\newtheorem{CorIntro}[ThmIntro]{Corollary}
\newtheorem{PropIntro}[ThmIntro]{Proposition}

\newtheorem{thm}{Theorem}[section]
\newtheorem{cor}[thm]{Corollary}
\newtheorem{lem}[thm]{Lemma}

\newtheorem{prop}[thm]{Proposition}
\theoremstyle{definition}
\newtheorem{defn}[thm]{Definition}

\theoremstyle{remark}
\newtheorem{rem}[thm]{Remark}

\numberwithin{equation}{section}

\newcommand{\Z}{\mathbf{Z}}

\newcommand{\N}{\mathbf{N}}
\newcommand{\R}{\mathbf{R}}

\newcommand{\Aut}{\text{Aut}}

\newcommand{\HH}{\mathcal{H}}

\newcommand{\bpr}{\noindent \textbf{Proof}: ~}

\newcommand{\epr}{~$\blacksquare$}

\newcommand{\Diam}{\textnormal{Diam}}

\title{Quantitative property~A, Poincaré inequalities, $L^p$-compression and $L^p$-distortion for metric measure spaces.}
\author{Romain Tessera}
\date{\today}

\begin{document}

\baselineskip=16pt

\maketitle

\begin{abstract}
We introduce a quantitative version of Property~A in order to
estimate the $L^p$-compressions of a metric measure space $X$. We
obtain various estimates for spaces with sub-exponential volume
growth. This quantitative property~A also appears to be useful to
yield upper bounds on the $L^p$-distortion of finite metric
spaces. Namely, we obtain new sharp results for finite subsets of
homogeneous Riemannian manifolds. We also introduce a general form
of Poincaré inequalities that provide constraints on compressions,
and lower bounds on distortion. These inequalities are used to
prove the optimality of some of our results.
\medskip

    \hfill\break \noindent {\sl Mathematics Subject
Classification:} Primary 51F99; Secondary 43A85. \hfill\break {\sl
Key words and Phrases:} Uniform embeddings, metric spaces, coarse
geometry, graphs, Banach spaces, volume growth, finite metric
spaces.
\end{abstract}

\pagebreak

\tableofcontents


\section{Introduction}

In \cite{Yu}, Yu introduced a {\it weak F\o lner} property for
metric spaces that he called Property~A. He proved that a metric
space satisfying this property uniformly embeds into a Hilbert
space. In \cite{Tu}, it is proved that a discrete metric space with
subexponential growth has Property~A and therefore, uniformly embeds
into a Hilbert space (here, we give a very short proof of this fact
when $X$ is assumed to be coarsely geodesic, e.g. if $X$ is a
graph). In this paper, we define a quantitative $L^p$-version of
Property~A and use it to obtain uniform embeddings of metric measure
spaces with subexponential growth into $L^p$ with compressions
satisfying some lower estimates.

Let us introduce some notation. The volume of the closed balls
$B(x,r)$ is denoted by $V(x,r)$. An $L^p$-space will mean a Banach
space of the form $L^p(\Omega,\mu)$ where $(\Omega,\mu)$ is some
measure space.

Let $f,g: \R_+\to\R_+$ be non-decreasing functions. We write
respectively $f\preceq g$, $f\prec g$ if there exists $C>0$ such
that $f(t)\leq Cg(Ct)+C$, resp. $f(t)=o(g(ct))$ for any $c>0$. We
write $f\approx g$ if both $f\preceq g$ and $g\preceq f$. The
asymptotic behavior of $f$ is its class modulo the equivalence
relation $\approx$.

The compression of a uniform embedding~\footnote{See the preliminary
section for a definition of a uniform embedding.} $F:X\to
L^p(\Omega,\mu)$ of a metric space $(X,d)$ into an $L^p$-space is
the following non-decreasing function $$\rho_F(t)=\inf_{d(x,y)\geq
t}\|F(x)-F(y)\|_p\quad \forall t>0.$$ A uniform embedding of a graph
into any other metric space is always large-scale Lipschitz, hence
$\rho_F(t)\leq Ct+C$ for some constant $C$. We are interested in
knowing how close the compression associated to a uniform embedding
can be from an affine function. Following \cite{GK}, for any $p\geq
1$, we can associate to a metric space $X$ a quasi-isometry
invariant quantity, denoted by $R_p$, by taking the supremum of all
positive $a$ such that there exists a uniform embedding $F:X\to
L^p(\Omega,\mu)$ with compression $\rho_F \geq ct^a.$

\subsection*{Quantitative Property~A and construction of uniform embeddings in $L^p$}

Let us give a definition\footnote{See also \cite[Definition
2.1]{R}.} of Yu's Property~A for metric measure spaces that
coincides with the usual one in the case of discrete metric spaces.

\begin{defn}\label{Adef}
We say that a metric measure space $X$ has Property~A if there
exists a sequence of families of probability densities on $X$:
$((\psi_{n,x})_{x\in X})_{n\in \N}$ such that
\begin{itemize}
\item[(i)] for every $n\in \N$, the support of each $\psi_{n,x}$ lies
in the (closed) ball $B(x,n)$ and \item[(ii)]
$\|\psi_{n,x}-\psi_{n,y}\|_1$ goes to zero when $n\to\infty$
uniformly on controlled sets $\{(x,y)\in X^2,\; d(x,y)\leq r\}$.
\end{itemize}
\end{defn}

The following proposition follows immediately from basic
$L^p$-calculus and its proof is left to the reader.

\begin{prop}\label{A_p}
Property~A is equivalent to the following statement. Let $1\leq
p<\infty$. There exists a sequence of families of unit vectors in
$L^p(X)$: $((\psi_{n,x})_{x\in X})_{n\in \N}$ such that
\begin{itemize}
\item[(i)] for every $n\in \N$, the support of each $\psi_{n,x}$ lies
in  $B(x,n)$ and \item[(ii)] $\|\psi_{n,x}-\psi_{n,y}\|_p$ goes to
zero when $n\to\infty$ uniformly on controlled sets $\{(x,y)\in
X^2,\; d(x,y)\leq r\}$.\epr
\end{itemize}
\end{prop}

The main conceptual tool in this paper the following quantitative
version in $L^p$ of Property A.

\begin{defn} Let $X=(X,d,\mu)$ be metric measure space, $J:\;\R_+\to\R_+$ be some increasing map and let $1\leq p<\infty$.
We say that $X$ has property $A(J,p)$ if for every $n\in \N$, there
exists a map $\psi_n:\;X\to L^p(X)$ such that
\begin{itemize}
\item for every $x\in X$, $\|\psi_{n,x}\|_p\geq J(n)$, \item
$\|\psi_{n,x}-\psi_{n,y}\|_p\leq d(x,y)$,
\item $\psi_{n,x}$ is supported in $B(x,n)$.
\end{itemize}
\end{defn}
The function $J$, that we call the A-profile in $L^p$, is a
increasing function dominated by $t$. This definition is motivated
by the following central observation.

\begin{PropIntro}\textnormal{(see Proposition~\ref{A/uniform})}
Let $X$ be a metric measure space satisfying Property $A(J,p)$.
Then, for every increasing function $f$ satisfying
\begin{equation}\tag{$J,p$}
\int_1^{\infty}\left(\frac {f(t)}{J(t)}\right)^p\frac{dt}{t}<\infty,
\end{equation}
there exists a large-scale Lipschitz uniform embedding $F$ of $X$
into $\oplus_{\infty}^{\ell^p} L^p(X,\mu)$ with compression
$\rho\succeq f$. In particular,
$$R_p(X)\geq \liminf_{t\to \infty}\log J(t).$$
\end{PropIntro}

We give estimates of the A-profile in $L^p$ for spaces with
subexponential growth.
\begin{PropIntro}\textnormal{(see Proposition~\ref{betterprop})}
Let $(X,d,\mu)$ be a metric measure space. Assume that there exists
an increasing function $v$ and constants $C\geq 1$ and $d>0$ such
that
$$1\leq v(r)\leq V(x,r)\leq Cv(r), \quad \forall x\in X, \forall r\geq 1;$$
Then for every $p\geq 1$, $X$ satisfies Equivariant Property
$A(J,p)$ with $J(t)\approx t/\log v(t)$.
\end{PropIntro}

For instance, we obtain the following corollary.
\begin{CorIntro}\textnormal{(see Proposition~\ref{bettercor})}
Keep the same hypothesis as in Proposition~\ref{betterprop} and
assume that $v(t)\preceq e^{t^{\beta}}$, for some $\beta<1$. Then,
for every $p\geq 1$,
$$R_p(X)\geq 1-\beta.~\blacksquare$$
\end{CorIntro}

Recall that a graph is called doubling if there exists a constant
$C$ such that $V(x,2r)\leq CV(x,r)$ for every $r>0$, $x\in X$. We
say that a graph is uniformly doubling if there exists an increasing
function $v$ satisfying $v(r)\leq V(x,r)\leq Cv(r)$ for every $r>0$,
$x\in X$, and a doubling property $v(2r)\leq C'v(r)$.

\begin{PropIntro}\textnormal{(see Proposition~\ref{doubling})}
Let $(X,d,\mu)$ be a uniformly doubling metric measure space. Then
for every $p\geq 1$, $X$ satisfies the property $A(J,p)$ with
$J(t)\approx t$.
\end{PropIntro}

\begin{CorIntro}\label{doubintro}\textnormal{(see Corollary~\ref{doublingcor})}
Let $X$ be uniformly doubling graph and let $p\geq 2$. Then, for
every increasing $f:\R_+\to\R_+$ satisfying
$$\int_1^{\infty}\left(\frac
{f(t)}{t}\right)^p\frac{dt}{t}<\infty,\quad (C_p)$$ there exists a
uniform embedding $F$ of $X$ into $\oplus_{\infty}^{\ell^p}
L^p(X,\mu)$ with compression $\rho_F\succeq f$. In particular,
$$R_p(X)=1.$$
\end{CorIntro}
\noindent{\bf About Condition $(C_p)$.} \begin{itemize}

\item First, note that if $ p\leq q$, then $(C_p)$ implies $(C_q)$:
this immediately follows from the fact that a nondecreasing function
$f$ satisfying $(C^p)$ also satisfies $f(t)/t=O(1).$

\item If $f$ and $h$ are two increasing functions such that
$f\preceq h$ and $h$ satisfies $(C_p)$, then $f$ satisfies $(C_p)$.

\item The function $f(t)=t^a$ satisfies $(C_p)$ for every $a<1$ but
not for $a=1$. More precisely, the function
$$f(t)=\frac{t}{(\log t)^{1/p}}$$
does not satisfy $(C_p)$ but
$$f(t)=\frac{t}{(\log t)^{a/p}}$$
satisfies $(C_p)$ for every $a>1$.

\item Let us call a function $f$ sublinear if $f(t)/t\to 0$ when
$t\to \infty$. Surprisingly, one can easily check \cite{compr} that
there exists no sublinear function that dominate all functions
satisfying Property $C_p$. Hence, by Corollary~\ref{doubintro}, a
function that dominates all the compression functions associated to
uniform embeddings of a uniformly doubling space into $L^p$ is at
least linear.
\end{itemize}

\medskip

In \cite{compr}, we proved that Corollary~\label{doubintro} is
actually true for a large variety of metric spaces, such as
homogeneous Riemannian manifolds, 3-regular trees, etc. In the
case of a 3-regular tree, the result is sharp since in turn, the
compression function associated to a uniform embedding in $L^p$
has to satisfy condition $(C_p)$. More surprising is that there
exists a doubling metric measure space having such a property.

\begin{PropIntro}\label{contrex}\textnormal{(see the remark preceding Proposition~\ref{contrexprop})}
There exists an infinite uniformly doubling graph such that for any
uniform embedding $F$ of $X$ into an Hilbert space, the compression
function of $F$ has to satisfy condition $(C_p)$, i.e.
$$\int_1^{\infty}\left(\frac
{\rho_F(t)}{t}\right)^2\frac{dt}{t}<\infty.$$
\end{PropIntro}

\begin{rem}

Note that Corollary~\ref{doubintro} should remain true if we merely
assume that $X$ is doubling as suggested by the result of
Assouad~\cite{A} that $R_p(X)=1$ for any doubling metric measure
space. On the other hand, this lack of generality is partly
compensated by the following ``equivariant" property of our
constructions.

\subsection*{Equivariance}

\begin{defn}\label{EquivEmbDef}
Let $(X,d)$ be a metric space. Consider a group $G$ of isometries of
$X$. We say that a map $F$ from $X$ to a metric space $Y$ is
$G$-equivariant if there exists an action of $G$ by isometries on
$Y$ such that $F$ commutes to theses actions, i.e. for all $x\in X$
and $g\in G,$
$$F(gx)=gF(x).$$
\end{defn}
All our constructions (so in particular in
Corollary~\ref{doubintro}) provide uniform embeddings that are
equivariant with respect to the group $Aut(X)$ of measure-preserving
isometries of the metric measure space $X$. In particular, if $G$ is
a group with subexponential growth, these constructions provide
proper isometric actions on $L^p$-spaces with the given compression
(see also \cite{compr}).
\end{rem}
\begin{rem}
Let us also emphasize the fact that our constructions are explicit
and involve relatively simple formulas.
\end{rem}

\subsection*{$L^p$-distorsion of finite metric spaces}

We also relate the quantitative property~A to the $L^p$-distortion
of finite metric spaces.

Our main result concerning the $L^p$-distorsion is that a finite
$1$-discrete\footnote{By $1$-discrete, we mean that the distance
between two points of $X$ is at least $1$.} subset of a homogeneous
manifold satisfies the following inequality
\begin{equation}\label{distupper}
c_p(X)\leq C (\log \Diam(X))^{1/p}
\end{equation}
where $p\geq 2$, and $C$ is a constant depending on the group. This
result is optimal in the sense that for any Lie group (or any
homogeneous Riemannian manifold) with exponential volume growth,
there exists an increasing sequence of $1$-discrete subsets $X_n$ of
diameter $n$ satisfying
\begin{equation}\label{distlower}
c_p(X_n)\geq c (\log \Diam(X_n))^{1/p}\quad \forall n\in \N.
\end{equation}
Note that Bourgain~\cite{Bourgain} proved the
inequalities~(\ref{distupper}) and (\ref{distlower}) in case $X$ is
a finite binary rooted tree. We deduce (\ref{distlower}) for Lie
groups with exponential growth from Bourgain's theorem and from the
fact~\cite{QI} that any Lie group with exponential growth contains a
quasi-isometrically embedded infinite binary rooted tree.

We also reprove~\cite[Theorem~4.1]{GKL} the optimal upper
bound~(\ref{distupper}) for the $L^p$-distortion of a
uniform\footnote{This bound actually holds for all doubling metric
spaces.} doubling metric spaces, the constant $C$ only depending
on the doubling constant of the metric space. Again, we loose some
generality by assuming uniform doubling property instead of
doubling property, but in counterpart, we get very explicit
embeddings, defined by simple expressions involving only the
metric\footnote{In \cite{GKL}, the constructions involve choices,
either of nets at various scales, for the Bourgain-style
embeddings, or of partitions for the Rao-style embeddings, which,
at least at first sight, prevent them from being equivariant.} and
the measure (in particular these embeddings are equivariant).

\subsection*{Optimality of the constructions and Poincar\'{e} inequalities.}

The graph of Proposition~\ref{contrex} is a planar self-similar
graph introduced and studied in \cite{L,L'}. In~\cite{GKL}, the
authors show that this graph satisfies a ``Poincaré-style"
inequality (for short, let us say Poincar\'e inequality) and they
deduce lower bounds on their $L^p$-distortions. Here, we use this
Poincaré inequality to prove Proposition~\ref{contrex}. The crucial
role of Poincaré inequalities for obtaining lower bounds on Hilbert
distortion has already been noticed in \cite{LMN}. Here, we try to
define the ``more general possible" Poincar\'e inequalities that
could be used to obtain, either constraints on the compression of
uniform embeddings into $L^p$-space, or lower bounds on the
$L^p$-distorsion, for $1\leq p<\infty$. We also propose a
generalization of these inequalities in order to treat uniform
embeddings into more general Banach spaces. We hope that these
definitions will be helpful in the future. In particular, proving
that Heisenberg satisfies a cumulated Poincaré
inequalities\footnote{See Section~\ref{Poincaresection}.} CP(J,p)
with $J(t)=ct$ would provide optimal constraints on its
$L^p$-compressions, for $p>1$. Another (weaker) consequence would be
that the $L^p$-distortion of balls of radius $r$ of the standard
Cayley graph of the discrete Heisenberg group is larger than $c(\log
r)^{\min(1/2,1/p)}$, which is not known, even for $p=2$, at least to
our knowledge.

\subsection*{Organization of the paper}

\begin{itemize}
\item In Section \ref{equivariantASection}, we introduce the
equivariant property A and give its interpretation in terms of the
quasi-regular representation of Aut(X) on $L^p(X)$. We also prove
that spaces with subexponential volume growth have equivariant
property~A.
\item The central part of the paper is
Section~\ref{GeomCondSection}. In
Section~\ref{Aquantitativesection}, we show how the A-profile can be
used to construct uniform embeddings with ``good" compression. In
Section~\ref{Poincaresection}, we introduce general forms of
Poincar\'{e} inequalities that provide constraints on the
compression of uniform embeddings.
\item Finally, in Section~\ref{ApplicationSection}, we
estimate the A-profile for spaces with subexponential volume growth
and Homogeneous Riemannian manifolds. Applying the results of
Section~\ref{GeomCondSection}, we obtain explicit constructions of
uniform embeddings of these spaces into $L^p$-spaces.
\end{itemize}

\section{Some preliminaries about uniform
embeddings}\label{preliminarysection}

In this section, we introduce the definitions of a uniform embedding
between metric spaces and of the compression function associated to
a uniform embedding.

Let $(X,d_X)$ and $(Y,d_Y)$ be metric spaces. A map $F: X\to Y$ is
called a uniform embedding of $X$ into $Y$ if there exists two
increasing, non bounded maps $\rho_1$ and $\rho_2$ such that
$$\rho_1(d_X(x,y))\leq d_Y(F(x),F(y))\leq \rho_2(d_X(x,y)).$$

A map $F: X\to Y$ is called a quasi-isometric embedding if $\rho_1$
and $\rho_2$ can be chosen affine (non-constant). The main purpose
of this paper is, given a metric space $X$, to find "good" uniform
embeddings of $X$ into some $L^p$-space. By good, we mean as close
as possible to a quasi-isometric embedding. Hence, the quality of a
uniform embedding will be measured by the asymptotics of $\rho_1$
and $\rho_2$. More precisely, let us define the compression of $F$
to be the supremum $\rho_F$ of all functions $\rho_1$ satisfying the
above inequality and the dilatation to be the infimum $\delta_F$ of
all functions $\rho_2$. Hence, $F$ is quasi-isometric if and only if
$\rho_F$  and $\delta_F$ are both asymptotically equivalent to
affine functions. To measure how far we are from this situation, one
can define the following function
$$\theta_F(t)=\exp \left( |\log(\rho_F(t)/t)|+
|\log(\delta_F(t)/t)|\right).$$ One can easily check that $F$ is a
quasi-isometric embedding if and only if $\theta_f$ is bounded.
Moreover, if $F$ is large-scale Lipschitz, i.e. if there is a
constant $C$ such that $\delta_F \leq Ct+C$, then $\theta_F(t)
\approx \rho_F(t)/t.$ The following well-known proposition shows
that this situation is actually very common. Recall that a metric
space $(X,d)$ is called {\it quasi-geodesic} if there exist $b>0$
and $\gamma\geq 1$ such that for all $x,y\in X$, there exists a
chain $x=x_0,x_1,\ldots,x_n=y$ satisfying
$$n\leq \gamma d(x,y),\textnormal{ and}$$
$$\forall k=1,\ldots,n,\quad d(x_{k-1},x_k)\leq b.$$
Such a chain is called an $b$-quasi-geodesic chain between $x$ and
$y$.

\begin{prop}
Let $X$ and $Y$ be two metric spaces such that $X$ is
quasi-geodesic. Then, any uniform embedding $F$ from $X$ to $Y$ is
large-scale Lipschitz.
\end{prop}
\bpr Let $x$ and $y$ be two elements of $X$, and let
$x=x_0,x_1,\ldots,x_n=y$ be an $b$-quasi-geodesic chain. Then,
$$d_Y(F(x),F(y))\leq (n+1)\delta(b)\leq \gamma\delta(b) (d(x,y)+1).~\blacksquare$$

\begin{rem}
In this paper, all the uniform embeddings that we construct are
large-scale Lipschitz, so we will focus on the compression function
$\rho_F$ instead of $\theta_F$.
\end{rem}

\begin{defn}\cite{GK}
Fix $p\geq 1$. The $L^p$-compression rate $R_p(X)$ of a metric space
$X$ is the supremum of $\alpha$ such that there exits a large-scale
Lipschitz uniform embedding from $X$ into a $L^p$-space with
compression $\rho(t)\succeq t^{\alpha}.$
\end{defn}

\begin{rem}
Note that $R_p(X)$ is invariant under quasi-isometry. More
generally, let $u:Y\to X$ be a quasi-isometric embedding from $X$ to
$Y$. Assume that $X$ admits a large-scale Lipschitz uniform
embedding $F$ into some $L^p$-space with compression $\rho_F$, then
$F\circ u$ defines a large-scale Lipschitz uniform embedding of $Y$
whose compression satisfies $\rho_{F\circ u}\succeq \rho_F$.
\end{rem}

\section{Equivariant Property~A}\label{equivariantASection}

\subsection{Equivariant Property~A and quasi-regular
representations of Aut$(X)$}

Let us denote by Aut$(X)$ the group of measure-preserving isometries
of $X$. We define a notion of ``equivariant" property~A, which means
that it behaves well under the action of Aut$(X)$. We will see that
this property implies that the quasi-regular representation of
Aut$(X)$ in $L^p(X)$ has almost invariant vectors for every $1\leq
p<\infty$.

\begin{defn}[\textnormal{Equivariant Property~A}]
Let $G$ be a group of isometries of $X$.  We say that a metric
measure space $X$ has $G$-equivariant Property~A if there exists a
sequence of families of unit vectors in $L^p(X)$ for one
(equivalently for any) $1\leq p<\infty$: $((\psi_{n,x})_{x\in
X})_{n\in \N}$ satisfying the conditions of Definition~\ref{Adef}
and the following additional one. For every $n\in \N$, $x,y\in X$
and $g\in G$,
\begin{equation}\label{equivar}
\psi_{n,g x}(y)=\psi_{n,x}(g^{-1}y).
\end{equation}
If $G$ is the entire group of isometries of $X$, then we just say
that $X$ has has Equivariant Property A (the same if $X=(X,d,\mu)$
is a metric measure space, and $G$ is the group of
measure-preserving isometries of $X$).
\end{defn}\label{equivAdef}

\begin{rem}
Note that if $f_{n,x}$ is defined only in terms of metric measure
properties around the point $x$, such as $V(x,r)$ or $1_{B(x,r)}$
where $r$ is a constant for instance, then it satisfies
(\ref{equivar}). This will be the case of all our constructions.
\end{rem}

\begin{prop}\label{equivprop}
Assume that $X$ has $G$-equivariant Property~A, then the
quasi-regular representation of $G$ on $L^p(X)$ for any $1\leq
p<\infty$ has almost invariant vectors. Moreover, if $G$ acts
transitively on $X$, then the converse is also true.
\end{prop}
\bpr Let us prove the first assertion for $p=1$. Let
$((\psi_{n,x})_{x\in X})_{n\in \N}$ satisfy the assumptions of
Definition~\ref{equivAdef}. Then by (\ref{equivar}), the sequence
$h_n=\psi_{n,x}$ for any fixed $x$ is almost-$G$-invariant.
Conversely, if $X$ is homogeneous and if $h_n$ is an
almost-$G$-invariant sequence in $L^p(X)$, then, given some $x_0\in
X$, we can define a sequence of families of unit vectors in
$L^p(X)$: $((\psi_{n,x})_{x\in X})_{n\in \N}$ satisfying the
conditions of Definition~\ref{equivAdef}, by
$\psi_{n,gx_0}(y)=h_n(g^{-1}y)$.\epr

\subsection{Equivariant Property~A for metric measure spaces with
subexponential growth}

In this section, we give a short proof of the fact that
subexponential growth implies Property~A. This is originally due to
Tu~\cite{Tu}. Tu's theorem works for any discrete metric space, so a
slight adaptation makes it work for any metric measure space. The
counterpart of this generality is that the proof is quite
complicated and does not yield any equivariance. Here, restricting
ourself to a certain class of metric measure spaces that includes
all graphs and Riemannian manifolds for instance, we give a short
proof that subexponential growth implies Equivariant~Property~A.

\

Recall that a metric measure space $X$ has bounded geometry if for
every $r>0$, there exists $C_r<\infty$ such that
$$C_r^{-1}\leq V(x,r)\leq C_r\quad \forall x\in X.$$
Consider some $b>0$ and define a $b$-geodesic distance on $X$ by
setting
$$d_b(x,y)=\inf_{\gamma}l(\gamma)$$
where $\gamma$ runs over all chains $x=x_0,\ldots, x_m=y$ such that
$d(x_{i-1},x_i)\leq b$ for all $1\leq i\leq m$, and where
$l(\gamma)=\sum_{i=1}^md(x_i,x_{i-1})$ denotes the length of
$\gamma$.

\begin{defn}
We say that a metric space $(X,d)$ is coarsely geodesic if there
exists $b>0$ such that the identity map $(X,d_b)\to (X,d)$ is a
uniform embedding.
\end{defn}

Let $(X,d,\mu)$ be a metric measure space such that $(X,d)$ is
coarsely geodesic. As $d\geq d_b,$ we have $V_b(x,r)\leq V(x,r)$,
where $V_b$ denotes the volume of balls in $(X,d_b,\mu).$

\begin{prop}\label{A}
Let $(X,d,\mu)$ be a coarsely geodesic metric measure space. Assume
that there exists a subexponential function $v:\; \R_+\to\R_+$ such
that $1\leq V(x,r)\leq v(r)$ for every $x\in X$ and $r\geq 0$. Then
$X$ has Equivariant Property~A.
\end{prop}

\bpr Since $X$ is coarsely geodesic, we can assume without loss of
generality that $d$ is a $1$-geodesic distance (replacing $d$ with
$d_1$). Denote $S_h(x,r)=V(x,r+h)-V(x,r).$ It is then easy to see by
a covering argument that for any $h>0$, there exists a constant
$C_h<\infty$ such that $V(x,r+h)\leq C_hV(x,r)$ for every $r>0$. We
define a sequence of families of probability densities
$(\psi_{n,x})$ by
$$\psi_{n,x}=\frac{1}{n}\sum_{k=1}^n\frac{1}{V(x,k)}1_{B(x,k)}\quad \forall x\in X.$$
Let $x$ and $y$ be such that $d(x,y)\leq h$, with $h\in \N^*$. We
have
\begin{eqnarray*}
\left\|\frac{1}{V(x,k)}1_{B(x,k)}-\frac{1}{V(y,k)}1_{B(y,k)}\right\|
& \leq &
\left\|\frac{1}{V(x,k)}(1_{B(x,k)}-1_{B(y,k)})\right\|\\
& &+V(y,k)\left|\frac{1}{V(x,k)}-\frac{1}{V(y,k)}\right|\\
& \leq & 2\frac{S_h(x,k)}{V(x,k)}\\
\end{eqnarray*}
Thus,
\begin{eqnarray*}
\| \psi_{n,x}-\psi_{n,y}\| & \leq &
\frac{1}{n}\sum_{k=1}^n\left\|\frac{1}{V(x,k)}1_{B(x,k)}-\frac{1}{V(y,k)}1_{B(y,k)}\right\|\\
& \leq & \frac{2}{n}\sum_{k=1}^n\frac{S_h(x,k)}{V(x,k)}\\
& \leq & \frac{2C_h}{n}\sum_{k=1}^n\frac{S_h(x,k)}{V(x,k+h)}\\
\end{eqnarray*}
But,
$$\frac{S_h(x,k)}{V(x,k+h)}=\sum_{i=0}^{h-1}\frac{S_1(x,k+i)}{V(x,k+h)}\leq \sum_{i=0}^{h-1}\frac{S_1(x,k+i)}{V(x,k+i+1)}.$$
Hence,
\begin{eqnarray*}
\| \psi_{n,x}-\psi_{n,y}\|  & \leq & \frac{2hC_h}{n}\sum_{k=1}^{n+h}\frac{S_1(x,k)}{V(x,k+1)}\\
&= & \frac{2hC_h}{n}\sum_{k=1}^{n+h}\frac{V(x,k+1)-V(x,k)}{V(x,k+1)}\\
& \leq & \frac{2hC_h}{n}\log
\left(\frac{V(x,n+h)}{V(x,1)}\right)\\
& \leq & \frac{2hC_h}{n}\log v(n+h)
\end{eqnarray*}
We conclude since $v$ is subexponential. \epr

\section{Geometric conditions to control compression and
distortion}\label{GeomCondSection}

\subsection{Quantitative Property~A, construction of uniform
embeddings and upper bounds on
distortion}\label{Aquantitativesection}

\begin{defn} Let $X=(X,d,\mu)$ be metric measure space, $J:\;\R_+\to\R_+$ be some increasing map and let $1\leq p<\infty$.
We say that $X$ has property $A(J,p)$ if for every $n\in \N$, there
exists a map $\psi_n:\;X\to L^p(X)$ such that
\begin{itemize}
\item for every $x\in X$, $\|\psi_{n,x}\|_p\geq J(n)$, \item
$\|\psi_{n,x}-\psi_{n,y}\|_p\leq d(x,y)$,
\item $\psi_{n,x}$ is supported in $B(x,n)$.
\end{itemize}
\end{defn}
\begin{rem}\label{1->p_rem}
Basic $L^p$-calculus shows that if $q\geq p\geq 1$, then Property
$A(J,q)$ implies Property $A(J,p)$ and Property $A(J,p)$ (only)
implies Property $A(J^{p/q},q)$.
\end{rem}
This definition is motivated by the following two propositions.

\begin{prop}\label{A/uniform}
Let $X$ be a metric measure space satisfying Property $A(J,p)$.
Then, for every increasing function $f$ satisfying
\begin{equation}\tag{$J,p$}
\int_1^{\infty}\left(\frac {f(t)}{J(t)}\right)^p\frac{dt}{t}<\infty,
\end{equation}
there exists a large-scale Lipschitz uniform embedding $F$ of $X$
into $\oplus_{\infty}^{\ell^p} L^p(X,\mu)$ with compression
$\rho\succeq f$. In particular,
$$R_p(X)\geq \liminf_{t\to \infty}\log J(t).$$
\end{prop}

\bpr Choose a sequence $(\psi_{n,x})$ like in Proposition \ref{A_p}.
Fix an element $o$ in $X$ and define
$$F(x)=\oplus^{\ell^p}_{k\in
\N}F_k(x)$$ where
$$F_k(x)=\left(\frac{f(2^k)}{J(2^k)}\right)(\psi_{2^k,x}-\psi_{2^k,o}).$$

The fact that $F$ exists and is Lipschitz follows from the fact that
Condition $(J,p)$ is equivalent to
$$\sum_{k}\left(\frac{f(2^k)}{J(2^k)}\right)^p<\infty.$$
Hence, a direct computation yields
\begin{eqnarray*}
\|F_k(x)-F_k(y)\|_p & \leq &
\left(\sum_{k}\left(\frac{f(2^k)}{J(2^k)}\right)^p\|\psi_{2^k,x}-\psi_{2^k,y}\|_p^p\right)^{1/p}\\
& \leq &
d(x,y)\left(\sum_{k}\left(\frac{f(2^k)}{J(2^k)}\right)^p\right)^{1/p}.
\end{eqnarray*}

On the other hand, since $\psi_{2^k,x}$ is supported in $B(x,2^k)$,
if $d(x,y)>2.2^k$, then the supports of $\psi_{2^k,x}$ and
$\psi_{2^k,y}$ are disjoint. Thus
$$\|F(x)-F(y)\|_p\geq \|F_k(x)-F_k(y)\|_p\geq\left(\|\psi_{2^k,x}\|_p^p+\|\psi_{2^k,y}\|_p^p\right)^{1/p}/J(2^k) f(2^{k})\geq 2^{1/p}f(2^k).$$
whenever $d(x,y)> 2.2^k$. So we are done.\epr

\begin{rem}
Note that it may happen that $\sum_k(1/J(2^k))^p=\infty$.
Nevertheless, as soon as $J$ is not bounded and $f(t)=o(J(t))$, one
can choose an increasing injection $i: \N\to \N$ such that
$$\sum_{n}\left(\frac{f(2^{i(n)})}{J(2^{i(n)})}\right)^p<\infty,$$
so that
$$F(x)=\oplus_{n}^{\ell^p}\left(\frac{f(2^{i(n)})}{J(2^{i(n)})}\right)(\psi_{2^{i(n)}}(x)-\psi_{2^{i(n)}}(o))$$
defines a uniform embedding of $X$ whose compression satisfies
$\rho\succeq f\circ i^{-1}$.
\end{rem}

\begin{prop}\label{A/distorsion}
Let $X$ be a finite metric space satisfying Property $A(J,p)$. Then,
$$c_p(X)\leq 2\left(\int_1^{\Diam(X)/4}\left(\frac{t}{J(t)}\right)^p\frac{dt}{t}\right)^{1/p}.$$
In particular, if $J(t)\geq ct$, then
$$c_p(X)\leq C\left(\log (\Diam(X))\right)^{1/p}.$$
\end{prop}

\bpr Fix an element $o$ in $X$, set $n=[\log(\Diam(X))/2]$ and
define
$$F(x)=\oplus^{\ell^p}_{k\in
\N}F_k(x)$$ where
$$F_k(x)=\left(\frac{2^k}{J(2^k)}\right)(\psi_{2^k,x}-\psi_{2^k,o}).$$
We have
\begin{eqnarray*}
\|F(x)-F(y)\|_p & \leq &d(x,y)\left(\sum_{k=0}^n
\left(\frac{2^k}{J(2^k)}\right)^p\right)^{1/p}\\
&\leq &
d(x,y)\left(\int_1^{\Diam(X)/2}\left(\frac{t}{J(t/2)}\right)^p\frac{dt}{t}\right)^{1/p}\\
&=&
2^{2/p}d(x,y)\left(\int_1^{\Diam(X)/4}\left(\frac{t}{J(t)}\right)^p\frac{dt}{t}\right)^{1/p}.
\end{eqnarray*}
On the other hand, since $\psi_{2^k,x}$ is supported in $B(x,2^k)$,
if $d(x,y)> 2.2^k$, then the supports of $\psi_{2^k,x}$ and
$\psi_{2^k,y}$ are disjoint. Thus
$$\|F(x)-F(y)\|_p\geq \|F_k(x)-F_k(y)\|_p\geq 2^{k}\left(\|\psi_{2^k,x}\|_p^p+\|\psi_{2^k,y}\|_p^p\right)^{1/p}/J(2^k) \geq 2^{1/p}2^k.$$
whenever $d(x,y)>2.2^k$. So we are done.\epr

\begin{rem}\label{equivariant}[{\bf Equivariance}]
This remark concerns the embeddings $F$ constructed in both
Propositions~\ref{A/uniform} and \ref{A/distorsion}.

Let $G=\Aut(X)$ be the group of measure preserving isometries of
$X$. Assume that in Propositions~\ref{A/uniform} and
\ref{A/distorsion}, the metric measure space $X=(X,d,\mu)$ actually
satisfies the Equivariant property A(J,p), i.e. if
$$\psi_{n,g x}(y)=\psi_{n,x}(g^{-1}y).$$ Then, the maps $F$ constructed in the proofs of those propositions are $G$-equivariant, according to
Definition~\ref{EquivEmbDef}. More precisely, there exists an affine
isometric action $\sigma_F$ of $G$ on $\oplus_{\infty}^{\ell^p}
L^p(X,\mu))$, whose linear part is the action by composition (which
is isometric since the elements of $G$ preserve the measure), such
that for every $g\in G$ and every $x\in X$,
$$\sigma_F(g) F(x)= F(gx).$$
In particular, Hence for every $g\in G$, we have
$$\forall x,y\in X, \quad \| F(gx)-F(gy)\|_p=\| F(x)-F(y)\|_p.$$

In particular, if $X=G$ is a compactly generated, locally compact
group, then $b(g)=F(g)-F(1)$ defines a $1$-cocycle of $G$ on the
infinite direct sum of the left regular representation (see
\cite{compr}).
\end{rem}

\subsection{Poincaré inequalities, constraints on uniform embeddings
and lower bounds on distortion}\label{Poincaresection}

In this section, we introduce general ``Poincar\'e-like"
inequalities in order to provide obstructions to embed a metric
space into an $L^p$-space, for $1\leq p<\infty$. The reader will
note that these inequalities are trivially inherited from a
subspace (that is a subset equipped with the induced metric).

In the sequel, let $(X,d)$ be a metric space and let $1\leq
p<\infty$. For any $r>0$, we denote
$$E_r=\{(x,y)\in X^2, \; d(x,y)\geq r\}.$$

\subsection*{The poincar\'e inequality P(J,p)}

\begin{defn}
Let $J:\; \R_+\to \R_+$ be a increasing function and let $r>0$. We
say that $X$ satisfies a Poincaré inequalities P(J,p) at scale $r$
if the following holds. There exists a Borel probability $P_r$ on
$E_r$ and a Borel probabilily $Q_r$ on $X^2$ such that for every
compactly supported continuous functions $\varphi:\; X\to \R$,
$$\int_{E_r}\left(\frac{|\varphi(x)-\varphi(y)|}{J(r)}\right)^pdP_r(x,y)\leq
\int_{E_1}\left(\frac{|\varphi(x)-\varphi(y)|}{d(x,y)}\right)^pdQ_r(x,y).$$
\end{defn}
This definition is motivated by the following simple proposition.

\begin{prop}\label{poincarecompr}
Let $p>1$. If a Metric space $X$ satisfies a Poincaré inequality
P(J,p) at scale $r$, then for any measurable large-scale Lipschitz
map $F$ of $X$ to a $L^p$-space, the compression $\rho_F$ has to
satisfy
$$\rho_F(t)\leq J(t)$$
for $t\leq r.$ If $X$ is finite, then,
$$c_p(X)\geq \frac{r}{J(r)}.$$
\end{prop}

\bpr Let $F:\; X\to L^p([0,1],\lambda)$ be a measurable large-scale
Lipschitz map. For almost every $t\in [0,1]$, the map
$F_t(x)=|F(x)(t)|^p$ defines a measurable map from $X$ to $\R$.
Applying the Poincaré inequality to this map and then integrating
over $t$ yields, by Fubini Theorem,
$$\int_{E_r}\left(\frac{\|F(x)-F(y)\|_p}{J(r)}\right)^pdP_r(x,y)\leq
\int_{E_1}\left(\frac{\|F(x)-F(y)\|_p}{d(x,y)}\right)^pdQ_r(x,y).$$
Now, the bounds for $\rho$ and $c_p(X)$ follow easily.\epr
\medskip

\noindent\textnormal{{\bf The skew cube inequality}.} In \cite{AGS},
upper bounds on the Hilbert compression rate are proved for a wide
variety of finitely generated groups including Thompson's group $F$,
$\Z\wr\Z$, etc. To show these bounds, they consider for all $n\in
\N$, injective group morphisms $j:\Z^n\to G$. Then, they focus on
the image, say $C_n$, of the $n$-dimensional cube $\{-1,1\}^n$. Let
$F$ be a uniform embedding from $G$ into a Hilbert space $\HH$. They
apply the well-known screw-cube inequality in Hilbert spaces to
$F(C_n)$. This inequality says that sum of squares of edges of a
cube is less or equal than the sum of squares of its diagonals. To
conclude something about the compression of $F$, they need an upper
bound (depending on $n$) on the length of diagonals of $C_n$ and a
lower bounds on the length of its edges. It is easy to check that
this actually remains to prove a Poincar\'e inequality P(J,2) for a
certain function $J$ (for instance, $J(t)=t^{1/2}\log t$ for
Thompson's group; and $J(t)=t^{3/4}$ for $Z\wr \Z$).

Let us briefly explain how one can deduce a Poincar\'e inequality
from the skew cube inequality. Denote by $\Delta_n$ the set of edges
of $C_n$ (seen as a cube embedded in $G$), and by $D_n$ the set of
diagonals. Assume that for all $(x,y)\in \Delta_n$, $d(x,y)\leq l_n$
and for all $(x,y)\in D_n$, $d(x,y)\geq L_n$, which actually means
that $D_n\subset E_{L_n}$. We have $|\Delta_n|=n2^{n-1}$ and
$|D_n|=2^{n-1}$. Take a function $\varphi: G\to \R$. The skew cube
inequality for the image of $C_n$ under $\varphi$ yields
$$\sum_{(x,y)\in D_n}|\varphi(x)-\varphi(y)|^2\leq \sum_{(x,y)\in \Delta_n}|\varphi(x)-\varphi(y)|^2.$$
An easy computation shows that this implies the following inequality
$$\frac{1}{|D_n|}\sum_{(x,y)\in D_n}\left(\frac{|\Delta_n|^{1/2}|\varphi(x)-\varphi(y)|}{l_n|D_n|^{1/2}}\right)^2\leq \sum_{(x,y)\in \Delta_n}\left(\frac{|\varphi(x)-\varphi(y)|}{d(x,y)}\right)^2,$$
which is nothing but P(J,2) with
$J(L_n)=l_n|D_n|^{1/2}/|\Delta_n|^{1/2}=l_nn^{1/2}.$

\medskip

\noindent{\bf Expanders.} Note that a metric space satisfying P(J,p)
with a constant function $J$ does not admit any uniform embedding in
any $L^p$-space. This is the case of families of expanders when
$p=2$. Recall that a sequence of finite graphs $(X_i)_{i\in \N}$ is
called a family of expanders if
\begin{itemize}
\item for every $i\in I$, the degree of $X_i$ is bounded by a
constant $d$;

\item the cardinal $|X_i|$ of $X_i$ tends to infinity when $n$
goes to infinity;

\item there is a constant $C>0$ such that for all $i\in I$, and
every function $f:X\to \R$,
$$\frac{1}{|X_i|^2}\sum_{(x,y)\in X_i^2}|f(x)-f(y)|^2\leq \frac{C}{|X_i|}\sum_{x~y}|f(x)-f(y)|^2.$$
\end{itemize}
As the volume of a ball of radius $r$ in $X_i$ is less than $d^r$,
for $i$ large enough, we have $$|E_r|\geq |X_i|(|X_i|-d^r)\geq
|X_i|^2/2.$$ Hence, the third property of expanders is equivalent to
property P(J,2), where $J=2C$, $P_r$ is the average over $E_r$, and
$Q_r$ is the average over $\{(x,y)\in X_i^2,\; d(x,y)=1\}\subset
E_1$.

\subsection*{The cumulated poincar\'e inequality CP(J,p)}

More subtle, the following definition will provide a finer control
on distortions and compressions.

\begin{defn}
Let $K:\; \R_+\to \R_+$ be a increasing function and let $r>0$. We
say that $X$ satisfies a cumulated Poincaré inequalities CP(J,p) at
scale $r$ if the following holds. There exist Borel probabilities
$P_{r,k}$ on $E_{2^k}$ for $k=1,2,\ldots [\log_2 r]$ and a Borel
probabilily $Q_r$ on $E_1$ such that for every measurable function,
$$\sum_{k=1}^{[\log_2 r]}\int_{E_{2^k}}\left(\frac{|\varphi(x)-\varphi(y)|}{J(2^k)}\right)^pdP_{r,k}(x,y)\leq
\int_{E_1}\left(\frac{|\varphi(x)-\varphi(y)|}{d(x,y)}\right)^pdQ_r(x,y).$$
\end{defn}

Here is the main application of these inequalities.

\begin{prop}\label{cumpoincarecompr}
Let $p>1$. If a Metric space $X$ satisfies a cumulated Poincaré
inequality CP(J,p) at scale $r$, then for any Lipschitz map $F$ of
$X$ to a $L^p$-space, the compression $\rho_F$ has to satisfy
$$\int_1^{r}\left(\frac{\rho_F(t)}{J(t)}\right)^{q}\frac{dt}{t}\leq 1,$$
where $q=\max(2,p)$. If $X$ is finite, then,
$$c_p(X)\geq \left(\int_1^r \left(\frac{t}{J(t)}\right)^p\frac{dt}{t}\right)^{\min(1/p,1/2)}.$$
\end{prop}
\bpr By a similar argument as for last proposition, we obtain the
following inequality
$$\sum_{k=1}^{[\log_2 r]}\int_{E_{2^k}}\left(\frac{\|F(x)-F(y)\|_p}{J(2^k)}\right)^pdP_{r,k}(x,y)\leq
\int_{E_1}\left(\frac{\|F(x)-F(y)\|_p}{d(x,y)}\right)^pdQ_r(x,y).$$
And, again, the proposition follows easily.\epr

\medskip

\noindent{\bf Trees and doubling graphs.} For a tree or for the
doubling graph of Proposition~\ref{contrex}, the fact that any
uniform embedding into an $L^p$-space satisfies Property $C_p$ is a
consequence of the inequality\footnote{proved in \cite{Bourgain} for
the trees, and in \cite{GKL} for the doubling graph.} CP(J,p) for
$J(t)=ct$.

Note that this inequality does not say anything for $p=1$ since the
$3$-regular tree $T$ admits a (trivial) bi-Lipschitz embedding into
$\ell^1$.

\subsection*{Poincar\'e inequalities with values in a Banach space}

To generalize these inequalities in order to treat embeddings into
more general Banach spaces, we can define P(J,p) (resp. CP(J,p))
with values in a Banach space $E$ to be the same inequalities
applied to elements in the Banach space $L^p(X,E)$ consisting of
functions $\varphi:\;X\to E$ such that $x\to \|\varphi(x)\|$ is in
$L^p(X)$. We equip $L^p(X,E)$ with the norm
$$\|\varphi\|_p=\left(\int_X\|\varphi(x)\|^pd\mu(x)\right)^{1/p}.$$

\medskip

In \cite{Bourgain}, Bourgain proves that if $E$ is a uniformly
$p$-convex Banach space, then the $E$-distorsion of the binary
tree $T_n$ of dept $n$ is more than a constant times $(\log
n)^{1/q}$ where $q=\max\{p,2\}$. To obtain this result, he
actually proves that $T_n$ satisfies CP(J,p), with values in $E$,
and with $J(t)\approx t$.

\medskip

In \cite{Laf}, Lafforgue constructs a sequence of expanders
satisfying P(J,2) with values in any uniformly convex Banach space
$E$, and with $J=constant$. In particular his expanders do not
uniformly embed into $E$.

\subsection*{Markov type inequalities}

A way to obtain Poincaré inequalities is to use Markov chains on
$X$. This idea was introduced by Ball \cite{Bal}, and was since
then used in various contexts (see for instance \cite{NPSS}). A
similar notion was introduced in \cite{LNP} under the name of
property of Markov convexity. The Markov convexity property can be
used to obtain cumulative Poincaré inequalities with $J(t)\approx
t$ (in particular, in \cite{LNP}, then use them for trees).

\section{Application to certain classes of metric
spaces}\label{ApplicationSection}

\subsection{Spaces with subexponential growth}
The proof of Proposition~\ref{A} yields the following proposition.
\begin{prop}\label{J1}
Let $(X,d,\mu)$ be a quasi-geodesic metric measure space with
bounded geometry. Assume that there exists a subexponential function
$v:\; \R_+\to\R_+$ such that $1\leq V(x,r)\leq v(r)$ for every $x\in
X$ and $r\geq 1$. Then $X$ has Equivariant Property $A(J_p,p)$ for
every $1\leq p<\infty$, with $J_p(t)\approx (t/\log v(t))^{1/p}$.
\end{prop}
\bpr The proof of Proposition~\ref{A} gives the result for $p=1$.
Then, by Remark~\ref{1->p_rem}, we deduce it for all $1\leq
p<\infty$.\epr

\begin{cor}\label{subexpcompressioncor}
Assume that $v(t)\preceq e^{t^{\beta}}$, for some $\beta<1$. Then,
for every $p\geq 1$, $$R_p(X)\geq (1-\beta)/p.~\blacksquare$$
\end{cor}

\

We can also improve Proposition \ref{J1} by assuming some uniformity
on the volume of balls.

\begin{prop}\label{betterprop}
Let $(X,d,\mu)$ be a metric measure space. Assume that there exists
an increasing function $v$ and constants $C\geq 1$ and $d>0$ such
that
$$1\leq v(r)\leq V(x,r)\leq Cv(r), \quad \forall x\in X, \forall r\geq 1;$$
Then for every $p\geq 1$, $X$ satisfies Equivariant Property
$A(J,p)$ with $J(t)\approx t/\log v(t)$.
\end{prop}
\begin{cor}\label{bettercor}
Keep the same hypothesis as in Proposition~\ref{betterprop} and
assume that $v(t)\preceq e^{t^{\beta}}$, for some $\beta<1$. Then,
for every $p\geq 1$,
$$R_p(X)\geq 1-\beta.~\blacksquare$$
\end{cor}

\noindent{\bf Proof of Proposition~\ref{betterprop}.} Define
$$k(n)=\sup\{k; \; v(n-k)\geq v(n)/2\}$$ and
$$j(n)=\sup_{1\leq j\leq n}{k(j)}.$$
We have
$$v(n)\geq 2^{n/j(n)}v(1)$$
which implies
$$j(n)\geq \frac{n}{\log v(n)}.$$
Let $q_n\leq n$ be such that $j(n)=k(q_n).$  Now define
$$\psi_{n,x}=\frac{1}{v(q_n)^{1/p}}\sum_{k=1}^{q_n-1}1_{B(x,k)}.$$
If $d(x,y)\leq 1$, we have
$$\|\psi_{n,x}-\psi_{n,y}\|_p^p\leq \frac{\mu(B(x,q_n))}{v(q_n)}\leq C.$$
On the other hand
$$\|\psi_{n,x}\|_p\geq j(n)\left(\frac{\mu(B(x,q_n-j(n)))}{v(q_n)}\right)^{1/p}\geq (1/2)^{1/p}\frac{n}{\log v(n)}$$
so we are done.\epr

\begin{prop}\label{doubling}
Let $(X,d,\mu)$ be a uniformly doubling metric measure space. Then
for every $p\geq 1$, $X$ satisfies the property $A(J,p)$ with
$J(t)\approx t$.
\end{prop}
\begin{cor}\label{doublingcor}
Let $(X,d,\mu)$ be a uniformly doubling metric measure space. Let
$p\geq 1$ and $q=\max(p,2).$ Then, for every increasing
$f:\R_+\to\R_+$ satisfying
$$\int_1^{\infty}\left(\frac
{f(t)}{t}\right)^q\frac{dt}{t}<\infty,$$ there exists a uniform
embedding $F$ of $X$ into $\oplus_{\infty}^{\ell^p} L^p(X,\mu)$ with
compression $\rho\succeq f$. In particular,
$$R_p(X)=1.~\blacksquare$$ Moreover, if $X$ is finite, then
$$c_p(X)\leq C\left(\log \Diam(X)\right)^{1/q}$$
where $C$ only depends on the doubling constant of $X$.
\end{cor}

\noindent{\bf Proof of Proposition~\ref{doubling}.} Define
$$\psi_{n,x}=\sum_{k=[n/2]}^{n}\frac{1}{v(k)^{1/p}}1_{B(x,k)}\quad \forall x\in X.$$
We have
$$\|\psi_{n,x}\|_p\geq \frac{n}{2}\left(\frac{V(x,n)}{V(x,[n/2])}\right)^{1/p}\geq \frac{n}{2C^{1/p}}$$
where $C$ is the doubling constant of $X$. On the other hand, for
$x$ and $y$ in $X$ at distance $d$, we want to prove that
\begin{eqnarray*}
\|\psi_{n,x}-\psi_{y,n}\|_p& \leq C"d
\end{eqnarray*}
for some constant $C"$. For obvious reasons, we can assume that
$d\leq n/4$. Hence, we can write this difference as follows
\begin{eqnarray*}
\psi_{n,x}-\psi_{y,n}& = & \sum_{i=1}^{2d}\left(\sum_{j=1}^{[n/4d]}
\frac{1}{v(2jd+i)^{1/p}}\left(1_{B(x,2jd+i)}-1_{B(y,2jd+i)}\right)\right).
\end{eqnarray*}
But note that for any $i,j$, $$B(x,2jd+i)\vartriangle
B(y,2jd+i)\subset B(x,2(j+1)d+i)\smallsetminus B(x,2jd+i).$$ Hence,
taking the absolute value, we obtain
\begin{eqnarray*}
|\psi_{n,x}-\psi_{y,n}|& \leq &
\sum_{i=1}^{2d}\left(\sum_{j=1}^{[n/4d]}
\frac{1}{v(2jd+i)^{1/p}}1_{B(x,2(j+1)d+i)\smallsetminus
B(x,2jd+i)}\right).
\end{eqnarray*}
Moreover, for distinct $j$ and $j'$, $B(x,2(j+1)+i)\smallsetminus
B(x,2jd+i)$ and $B(x,2(j'+1)d+i)\smallsetminus B(x,2j'd+i)$ are
disjoint. Hence, taking the $L^p$-norm, we get
\begin{eqnarray*}
\|\psi_{n,x}-\psi_{y,n}\|_p& \leq &
\sum_{i=1}^{2d}\left(\sum_{j=1}^{[n/4d]}
\frac{1}{v(jd+i)}|B(x,2(j+1)d+i)\smallsetminus B(x,2jd+i)|\right)^{1/p}\\
&\leq& \sum_{i=1}^{2d}
\frac{1}{v([n/4]^{1/p})}\left(\sum_{j=1}^{[n/4d]}|B(x,2(j+1)d+i)\smallsetminus B(x,2jd+i)|\right)^{1/p}\\
&=& d\frac{V(x,2n+2d)^{1/p}}{v([n/2]^{1/p})}\\
&=& d\frac{V(x,2n)^{1/p}}{v([n/2]^{1/p})}\\
&=& C"d.~\blacksquare
\end{eqnarray*}

Now, the optimality of Corollary~\ref{doublingcor} follows from
Proposition~\ref{cumpoincarecompr} and  the following result,
essentially remarked in \cite[Theorem~4.1]{GKL}.

\begin{prop}\label{contrexprop}
There exists an infinite uniformly doubling graph satisfying CP(J,2)
for a linear increasing function $J$.
\end{prop}

\subsection{Homogeneous Riemannian manifolds}

\begin{thm}
Let $X$ be a homogeneous Riemannian manifold. Then, $X$ has Property
A(J,p) for all $p\geq 1$ and $J=ct$ where $c$ depends only on $X$.
Moreover, $X$ has Equivariant Property A if and only if Isom$(X)$ is
amenable.
\end{thm}
\bpr In \cite{compr}, we prove that for any amenable Lie group $G$,
equipped with a left Haar measure and with a left-invariant
Riemannian metric\footnote{Actually, we prove it using a word length
metric on $G$, but such a metric is quasi-isometric to any
Riemannian one.} and for every $1\leq p<\infty$ and every $n\in \N$,
there exists a measurable function $h_n:\; G\to \R$ whose support
lies in the unit ball of radius $n$ and such that for every element
$g\in G$ of length less than $1$,
$$\|h_n(g\cdot)-h_n\|_p\leq 1$$
and
$$\|h_n\|_p\geq cn$$
for a constant $c$ only depending on $G$. To see that $G$ satisfies
equivariant A(J,p) with $J(t)=ct$, we construct as in the proof of
Proposition~\ref{equivprop}, a sequence
$$\psi_{n,g}(y)=h_n(g^{-1}y).$$
Now, let $X$ be a homogeneous manifold and let $G$ be its group of
isometries. We have $X=G/K$ where $K$ is a compact subgroup of $G$.

First, assume that $G$ is amenable. Averaging them over $K$, we can
assume that the $h_n$ are $K$-bi-invariant and then we can push them
through the projection $G\to X$. We therefore get Equivariant
Property A(J,p) for $X$. Conversely, if $X$ has Equivariant Property
A, then by Proposition~\ref{equivprop}, the quasi-regular
representation of $G$ on $L^p(X)$ has almost-invariant vectors.
Lifting them to $G$, we obtain almost-invariant vectors on the left
regular representation of $G$ on $L^p(G)$, which implies that $G$ is
amenable.

We are left to prove that even when $G$ is not amenable, $X$
satisfies A(J,p). Recall that every connected Lie group has a
connected solvable co-compact subgroup. So any homogeneous
Riemannian manifold is actually quasi-isometric to some amenable
connected Lie group. Hence the first statement of the theorem
follows from Lemma~\ref{lemQI}.\epr
\begin{lem}\label{lemQI}
If $F: \; X\to Y$ is a quasi-isometry between two metric measure
spaces with bounded geometry, then if $X$ satisfies A(J,p), then $Y$
satisfies A(J',p) with $J'\approx J$ for some constant $c>0$.
\end{lem}
\noindent{\bf Proof the lemma.} If we assume that $X$ and $Y$ are
discrete spaces equipped with the counting measure and that $F$ is a
bi-Lipschitz map, then the claim is obvious. Now, to reduce to this
case, we just have to prove that we can replace $X$ and $Y$ by any
of their nets, which is essentially proved in
\cite[Lemma~2.2]{R}.\epr

\bigskip
\footnotesize

\noindent \noindent Romain Tessera\\
Department of mathematics, Vanderbilt University,\\ Stevenson
Center, Nashville, TN 37240 United,\\ E-mail:
\url{tessera@clipper.ens.fr}


\begin{thebibliography}{KM98b}


\bibitem[A]{A} P. {\sc Assouad}.
\newblock {\em Plongements lipschitziens dans $\R^n$}. \newblock Bull. Soc. Math. France,
111(4), 429-448, 1983.



\bibitem[AGS]{AGS} G. N. {\sc Arzhantseva}, V.S. {\sc Guba} and  M.
V. {\sc Sapir}. \newblock {\em Metrics on diagram groups and uniform
embeddings in Hilbert space}. ArXiv GR/0411605, 2005.


\bibitem[Ba]{Bal} K. {\sc Ball}.
\newblock {\em Markov chains, Riesz transforms and Lipschitz maps}. \newblock
Geom. Funct. Anal., 2 (2), 137-172, 1992.



\bibitem[Bou]{Bourgain} J. {\sc Bourgain}.
\newblock {\em The metrical interpretation of superreflexivity in Banach spaces}. \newblock Israel J. Math. ,
Vol. 56, N°2, 1986.

\bibitem[CT]{QI} Y. {\sc de Cornulier} and R. {\sc Tessera}.
\newblock {\em Quasi-isometrically embedded trees}. \newblock In preparation, 2005.

\bibitem[GK]{GK} E. {\sc Guentner}, J. {\sc Kaminker}.
\newblock {\em Exactness and uniform embeddability of discrete
groups}. \newblock J. London Math. Soc. {\bf 70}, 703-718, 2004.


\bibitem[GKL]{GKL} A. {\sc Gupta}, R. {\sc Krauthgamer} and J.R. {\sc Lee}.
\newblock {\em Bounded geometry, fractals, and low-distortion embeddings}. \newblock Proc. of the 44th Annual IEEE Symposium on Foundations of Computer Science, 2003.

\bibitem[L]{L} T.J. {\sc Laakso}. \newblock {\em Ahlfors $Q$-regular spaces with arbitrary $Q>1$ admiting weak Poincaré
inequality.}
\newblock Geom. Funct. Ann. 10(1), 111-123, 2000.

\bibitem[L']{L'} T.J. {\sc Laakso}. \newblock {\em Plane with $A_{\infty}$-weighted metric not bi-Lischitz embeddable to $\R^n$.}
\newblock Bull. London Math. Soc. 34(6), 667-676, 2002.

\bibitem[La]{Laf} V. {\sc Lafforgue}. \newblock {\em Un renforcement de la propriété (T).}
\newblock Preprint, 2006.


\bibitem[NPSS]{NPSS}  A. {\sc Naor}, Y. {\sc Perez}, O. {\sc Schramm}, S. {\sc Sheffield}.
\newblock {\em Markov chains in smooth Banach spaces and Gromov-hyperbolic metric spaces}. \newblock 
Duke Math. J., 134 (1), 165-197, 2006.


\bibitem[LMN]{LMN} N. {\sc Linial}, A. {\sc Magen} and A. {\sc Naor}. \newblock {\em Girth and Euclidean distortion.}
\newblock Proc. of the thiry-fourth annual ACM symposium on Theory of computing, 705-711, 2002.


\bibitem[LNP]{LNP} J.R. {\sc Lee}, A. {\sc Naor}, Y. {\sc Perez}.
\newblock {\em Trees and Markov convexity}. \newblock To appear in
Geom. Funct. Anal.


\bibitem[R]{R} J. {\sc Roe}. \newblock {\em Warped cones and Property~A.} \newblock Geom. and Topol. Pub.  {\bf 9}, 163-178, 2005.

\bibitem[T]{compr} R. {\sc Tessera}. \newblock {\em Asymptotic isoperimetry on groups and uniform embeddings into Banach spaces.} \newblock math.GR/0603138, 2006.


\bibitem[Tu]{Tu} J.L. {\sc Tu}. \newblock {\em Remarks on Yu's
``Property~A" for discrete metric spaces and groups.} \newblock
Bull. Soc. Math. France {\bf 129}, 115-139, 2001.


\bibitem[Yu]{Yu} G. {\sc Yu}. \newblock {\em The coarse Baum-Connes conjecture for spaces which admit a uniform embedding into Hilbert space.}
\newblock Invent. Math. {\bf 139}, 201-240, 2000.


\end{thebibliography}
\end{document}